\documentclass[a4paper,11pt]{article}
\usepackage{amsmath}
\usepackage{amsfonts}
\usepackage{amssymb}
\usepackage{amsthm}
\usepackage{hyperref}

\title{Combinatorics and topology \\of toric arrangements defined by root systems}
\author{Luca Moci}

\begin{document}

\newtheorem{te}{Theorem}
\newtheorem{lem}[te]{Lemma}
\newtheorem{co}[te]{Corollary}
\theoremstyle{definition}
\newtheorem{re}[te]{Remark}

\newcommand{\faktor}{\frac}
\newcommand{\PT}{\mathcal{C}_{0}(\Phi)}
\newcommand{\Pt}{\mathcal{C}_{0}(\Theta)}
\newcommand{\Tquo}{{\mathfrak{h}}/{\langle\Phi^\vee\rangle}}
\newcommand{\sym}{\mathfrak{S}}
\newcommand{\rp}{\mathcal{R}_\Phi}

\maketitle

%

\begin{abstract}
Given the toric (or toral) arrangement defined by a root system $\Phi$, we describe the poset of its \emph{layers} (connected components of intersections) and we count its elements. Indeed we show how to reduce to 0-dimensional layers, and in this case we provide an explicit formula involving the maximal subdiagrams of the affine Dynkin diagram of $\Phi$.
Then we compute the Euler characteristic and the Poincar\'{e} polynomial of the complement of the arrangement, which is the set of regular points of the torus.\\
\end{abstract}

\section{Introduction}
Let $\mathfrak{g}$ be a semisimple Lie algebra of rank $n$ over
$\mathbb{C}$, $\mathfrak{h}$ a Cartan subalgebra,
$\Phi\subset\mathfrak{h^*}$ and $\Phi^\vee\subset\mathfrak{h}$
respectively the root and coroot systems. The equations
$\{\alpha(h)=0,\:\alpha\in\Phi\}$ define in $\mathfrak{h}$ a
family $\mathcal{H}$ of intersecting hyperplanes. Let
$\langle\Phi^\vee\rangle$ be the lattice spanned by the
coroots: the quotient $T\doteq \Tquo$ is a complex torus of
rank $n$. Each root $\alpha$ takes integer values on
$\langle\Phi^\vee\rangle$, hence it induces a map $T\rightarrow
{\mathbb{C}}/{\mathbb{Z}}\simeq \mathbb{C^*}$ that we denote by
$e^{\alpha}$. This is a character of $T$; let $H_\alpha$ be its
kernel:
$$H_\alpha\doteq\{t\in T \:|\:e^\alpha(t)=1\}.$$

In this way $\Phi$ defines in $T$ a finite family of
hypersurfaces
$$\mathcal{T}\doteq\{H_\alpha,\:\alpha\in\Phi^+\}$$
(since clearly $H_\alpha=H_{-\alpha}$).
 $\mathcal{H}$ and $\mathcal{T}$ are called respectively the \emph{hyperplane arrangement} and the \emph{toric arrangement} defined by $\Phi$ (see for instance \cite{DPt}, \cite{li}, \cite{OT}). We call \emph{spaces} of $\mathcal{H}$ the intersections of elements of $\mathcal{H}$, and \emph{layers} of $\mathcal{T}$ the connected components of the intersections of elements of $\mathcal{T}$. We denote by $\mathcal{L}(\Phi)$ the set of the spaces of $\mathcal{H}$, by $\mathcal{C}(\Phi)$ the set of the layers of $\mathcal{T}$, and by $\mathcal{L}_d(\Phi)$ and $\mathcal{C}_d(\Phi)$ the sets of
$d-$dimensional spaces and layers. Clearly if
$\Phi=\Phi_1\times\Phi_2$ then
$\mathcal{L}(\Phi)=\mathcal{L}(\Phi_1)\times
\mathcal{L}(\Phi_2)$ and
$\mathcal{C}(\Phi)=\mathcal{C}(\Phi_1)\times
\mathcal{C}(\Phi_2)$, hence from now on we will suppose $\Phi$
to be irreducible.
Let $W$ be the Weyl group of $\Phi$: since $W$ permutes the roots, its natural action on $T$ restricts to an action on $\mathcal{C}(\Phi)$.\\

$\mathcal{H}$ is a classical object, whereas $\mathcal{T}$ has
recently been shown (\cite{DPt}) to provide a geometric way to
compute the values of the Kostant partition function. This
function counts in how many ways an element of the lattice
$\langle \Phi \rangle$ can be written as sum of positive roots,
and plays an important role in representation theory, since (by
Kostant's and Steinberg's formulae \cite{Ko}, \cite{St}) it
yields efficient computation of weight multiplicities and
Littlewood-Richardson coefficients, as shown in \cite{Co} using
results from \cite{c1}, \cite{c2}, \cite{c3}, \cite{c4}. The
values of Kostant partition function can be computed as a sum
of contributions given by the elements of $\PT$ (see \cite[Teor
3.2]{Co}).

Furthermore, let $\rp$ be the complement in $T$ of the union of all elements of $\mathcal{T}$. $\rp$ is known as the set of the \emph{regular points} of the torus $T$ and has been widely studied (see in particular \cite{DPt}, \cite{L1}, \cite{L2}). The cohomology of $\rp$ is direct sum of contributions given by the elements of $\mathcal{C}(\Phi)$ (see for instance \cite{DPt}). Then by describing the action of $W$ on $\mathcal{C}(\Phi)$ we implicitly obtain a $W-$equivariant decomposition of the cohomology of $\rp$, and by counting and classifying the elements of $\mathcal{C}(\Phi)$ we can compute the Poincar\'{e} polynomial of $\rp$.\\

We say that a subset $\Theta$ of $\Phi$ is a \emph{subsystem}
if it satisfies the following conditions:
\begin{enumerate}
  \item $\alpha\in\Theta\Rightarrow -\alpha\in\Theta$
  \item $\alpha,\beta\in\Theta$ and $\alpha+\beta\in\Phi
      \Rightarrow \alpha+\beta\in\Theta.$
\end{enumerate}
For each $t\in T$ let us define the following subsystem of
$\Phi$:
$$\Phi(t)\doteq\{\alpha\in\Phi|e^\alpha(t)=1\}.$$
and denote by $W(t)$ the stabilizer of $t$.

The aim of Section 2 is to describe $\PT$, which is the set of
points $t\in T$  such that $\Phi(t)$ has rank $n$. We call its
elements the \emph{points} of the arrangement $\mathcal{T}$.
Let $\alpha_1,\dots,\alpha_n$ be simple roots of $\Phi$,
$\alpha_0$ the lowest root (i.e. the opposite of the highest
root), and $\Phi_p$ the subsystem of $\Phi$ generated by
$\{\alpha_i\}_{0\leq i\leq n, i\neq p}$. Let $\Gamma$ be the
affine Dynkin diagram of $\Phi$ and $V(\Gamma)$ the set of its
vertices (a list of such diagrams can be found for instance in
\cite{He} or in \cite{KL}). $V(\Gamma)$ is in bijection with
$\{\alpha_0,\alpha_1,\dots,\alpha_n\}$, hence we can identify
each vertex $p$ with an integer from 0 to $n$. The diagram
$\Gamma_p$ obtained by removing from $\Gamma$ the vertex $p$
(and all adjacent edges) is the ordinary Dynkin diagram of
$\Phi_p$.
  Let $W_p$ be the Weyl group of $\Phi_p$, i.e. the subgroup of $W$ generated by all the reflections $s_{\alpha_0},\dots, s_{\alpha_n}$ except $s_{\alpha_p}$. Notice that $\Gamma_0$ is the Dynkin diagram of $\Phi$ and $W_0=W$.

Then we prove:
\begin{te}\label{pts}
There is a bijection between the $W-$orbits of $\PT$ and the
vertices of $\Gamma$, having the property that for every point
$t$ in the orbit $\mathcal{O}_p$ corresponding to the vertex
$p$, $\Phi(t)$ is $W-$conjugate to $\Phi_p$ and $W(t)$ is
$W-$conjugate to $W_p$.
\end{te}
As a corollary we get the formula
\begin{equation}
|\PT|=\sum_{p\in V(\Gamma)}\cfrac{|W|}{|W_p|}.\\
\end{equation}

In Section 3 we deal with layers of arbitrary dimension. For
each layer $C$ of $\mathcal{T}$ we consider the subsystem of
$\Phi$ $$\Phi_C\doteq \{\alpha\in\Phi | e^{\alpha}(t)=1 \;
\forall t\in C\}$$ and its \emph{completion}
$\overline{\Phi_C}\doteq \langle\Phi_C\rangle_{\mathbb{R}}\cap
\Phi$.

Let $\mathcal{K}_d$ be the set of subsystems $\Theta$ of $\Phi$
of rank $n-d$ that are \emph{complete} (i.e. such that
$\Theta=\overline{\Theta}$), and let
$\mathcal{C}^{\Phi}_{\Theta}$ be the set of layers $C$ such
that $\overline{\Phi_C}=~\Theta$. This gives a partition of the
layers:
$$\mathcal{C}_{d}(\Phi)=\bigsqcup_{\Theta\in \mathcal{K}_d} \mathcal{C}^{\Phi}_{\Theta}.$$
Notice that the subsystem of roots vanishing on a space of
$\mathcal{H}$ is always complete; then $\mathcal{K}_d$ is in
bijection with $\mathcal{L}_d$. The elements of $\mathcal{L}_d$
are classified and counted in \cite{OS}, \cite{OT}. Thus the
description of the sets $\mathcal{C}^{\Phi}_{\Theta}$ given in
Theorem \ref{cft} yields a classification of the layers of
$\mathcal{T}$. In particular we show that
$|\mathcal{C}^{\Phi}_{\Theta}|=n_\Theta^{-1}
|\mathcal{C}_{0}(\Theta)|$, where $n_\Theta$ is a natural
number depending only on the conjugacy class of $\Theta$, and
then
\begin{equation*}
|\mathcal{C}_{d}(\Phi)|=\sum_{\Theta\in \mathcal{K}_d} n_\Theta^{-1} |\mathcal{C}_{0}(\Theta)|.
\end{equation*}

In Section 4, using results of \cite{DPt} and \cite{ci}, we
deduce from Theorem \ref{pts} that the Euler characteristic of
$\rp$ is equal to $(-1)^n |W|$. Moreover, Corollary \ref{nla}
yields a formula for the Poincar\'{e} polynomial of $\rp$:
\begin{equation*}
P_\Phi(q)=\sum_{d=0}^n (-1)^d(q+1)^d q^{n-d }\sum_{\Theta\in \mathcal{K}_d}n_\Theta^{-1}|W^\Theta|.
\end{equation*}
By this formula $P_\Phi(q)$ can be explicitly computed.\\

I would like to thank Gus Lehrer and Paolo Papi for giving me useful suggestions, and Filippo Callegaro, Francesca Mori and Alessandro Pucci for some stimulating discussions. I also want express my gratitude to my supervisor Corrado De Concini for suggesting to me many key ideas, encouraging me, and giving me a lot of helpful advice.

\section{Points of the arrangement}
\subsection{Statements}
For all facts about Lie algebras and root systems we refer to
\cite{HA}. Let
$$\mathfrak{g}=\mathfrak{h}\oplus\bigoplus_{\alpha\in\Phi}\mathfrak{g}_\alpha$$
be the Cartan decomposition of $\mathfrak{g}$, and let us
choose nonzero elements $$X_0,X_1,\dots,X_n$$ in the
one-dimensional subalgebras
$\mathfrak{g}_{\alpha_0},\mathfrak{g}_{\alpha_1},\dots,\mathfrak{g}_{\alpha_n}$:
since $[\mathfrak{g}_\alpha,\mathfrak{g}_{\alpha
'}]=\mathfrak{g}_{\alpha+\alpha '}$ whenever $\alpha,\alpha
',\alpha+\alpha '\in\Phi$, we have that $X_0,X_1,\dots,X_n$
generate $\mathfrak{g}$. Let $a_0=1$ and for $p=1,\dots,n$ let
$a_p$ be the coefficient of $\alpha_p$ in $-\alpha_0$. For each
$p=0,\dots,n$ we define an automorphism $\sigma_p$ of
$\mathfrak{g}$ by
$$\sigma_p(X_j)\doteq \begin{cases}
X_j & \mbox{ if } j\neq p\\
e^{2\pi i  a_p^{-1}} X_j & \mbox{ if } j=p
\end{cases}$$

Let $G$ be the semisimple and simply connected linear algebraic
group having root system $\Phi$; then $\mathfrak{g}$ is the Lie
algebra of $G$, and $T$ is the maximal torus of $G$
corresponding to $\mathfrak{h}$ (see for instance \cite{Hu}).
$G$ acts on itself by conjugacy, and for each $g\in G$ the map
$k\mapsto gkg^{-1}$ is an automorphism of $G$. Its differential
$Ad(g)$ is an automorphism of $\mathfrak{g}$.

\begin{re}\label{re2}
For every $t\in\mathcal{C}_0(\Phi)$, let $\mathfrak{g}^{Ad(t)}$
be the subalgebra of the elements fixed by $Ad(t)$. For every
$\alpha\in\Phi$ and for every
$X_\alpha\in\mathfrak{g}_{\alpha}$ we have that
$$Ad(t)(X_\alpha)=~e^\alpha(t)X_\alpha$$
and then
$$\mathfrak{g}^{Ad(t)}=\mathfrak{h}\oplus\bigoplus_{\alpha\in\Phi(t)}\mathfrak{g}_\alpha.$$
On the other hand $\mathfrak{g}^{\sigma_p}$ is generated by the
subalgebras $\{\mathfrak{g}_{\alpha_i}\}_{0\leq i\leq n, i\neq
p}$. Then $\mathfrak{g}^{Ad(t)}$ and $\mathfrak{g}^{\sigma_p}$
are semisimple algebras having root system respectively
$\Phi(t)$ and $\Phi_p$. Our strategy will be to prove that for
each $t\in\mathcal{C}_0(\Phi)$, $Ad(t)$ is conjugate to some
$\sigma_p$. This implies that $\mathfrak{g}^{Ad(t)}$ is
conjugate to $\mathfrak{g}^{\sigma_p}$ and then $\Phi(t)$ to
$\Phi_p$, as claimed in Theorem \ref{pts}.
\end{re}
Then we want to give a bijection between vertices of $\Gamma$
and $W-$orbits of $\PT$ showing that, for every $t$ in the
orbit $\mathcal{O}_p$,  $Ad(t)$ is conjugate to $\sigma_p$.
However, since some of the $\sigma_p$ (as well as the
corresponding $\Phi_p)$ are themselves conjugate, this
bijection is not canonical. To make it canonical we should
merge the orbits corresponding to conjugate automorphisms: for
this we consider the action of a larger group.

Let $\Lambda(\Phi)\subset\mathfrak{h}$ be the lattice of the
\emph{coweights} of $\Phi$, i.e.
$$\Lambda(\Phi)\doteq\{h\in \mathfrak{h} | \alpha (h)\in\mathbb{Z}\;\forall \alpha\in \Phi\}.$$
The lattice spanned by the coroots $\langle\Phi^\vee\rangle$ is
a sublattice of $\Lambda(\Phi)$; set $$Z(\Phi)\doteq
\frac{\Lambda(\Phi)}{\langle\Phi^\vee\rangle}.$$ This finite
subgroup of $T$ coincides with $Z(G)$, the \emph{center} of
$G$. It is well known (see for instance \cite[13.4]{Hu}) that
\begin{equation}\label{adg}
    Ad(g)=id_{\mathfrak{g}} \Leftrightarrow g\in Z(\Phi).
\end{equation}

Notice that
$$Z(\Phi)=\{t\in T|\Phi(t)=\Phi\}$$
thus $Z(\Phi)\subseteq \PT$. Moreover, for each $z\in Z(\Phi),
t\in T, \alpha\in\Phi$,
$$e^\alpha(zt)=e^\alpha(z) e^\alpha(t)=e^\alpha(t)$$
and therefore $\Phi(zt)=\Phi(t)$. In particular $Z(\Phi)$ acts
by multiplication on $\PT$.
Notice that this action commutes with that of $W$: indeed, let
$$N\doteq N_G(T)$$ be the normalizer of $T$ in $G$. We recall
that $W\simeq N/T$ and the action of $W$ on $T$ is induced by
the conjugacy action of $N$. The elements of $Z(\Phi)=Z(G)$
commute with the elements of $G$, hence in particular with the
elements of $N$. Thus we get an action of $W\times Z(\Phi)$ on
$\PT$.

Let $Q$ be the set of the $Aut(\Gamma)$-orbits of $V(\Gamma)$.
If $p, p'\in V(\Gamma)$ are two representatives of $q\in Q$,
then $\Gamma_{p}\simeq\Gamma_{p'}$, thus $W_p\simeq W_{p'}$.
Moreover we will see (Corollary \ref{kco}(ii)) that $\sigma_p$
is conjugate to $\sigma_{p'}$. Then we can restate Theorem
\ref{pts} as follows.

\begin{te}\label{pca}~
There is a canonical bijection between $Q$ and the set of
$W\times Z(\Phi)-$orbits in $\PT$, having the property that if
$p\in V(\Gamma)$ is a representative of $q\in Q$, then:
\begin{enumerate}
\item every point $t$ in the corresponding orbit
    $\mathcal{O}_q$ induces an automorphism conjugate to
    $\sigma_p$;
\item the stabilizer of $t\in\mathcal{O}_q$ is isomorphic
    to $W_p\times Stab_{Aut(\Gamma)}p$.
\end{enumerate}
\end{te}

This theorem implies immediately the formula:
\begin{equation}|\PT|=\sum_{q\in Q}|q|\:\cfrac{|W|}{|W_p|}
\end{equation}
where $p$ is any representative of $q$. This is clearly equivalent to formula (1).\\

\begin{re}
If we view the elements of $\Lambda(\Phi)$ as translations, we
can define a group of isometries of $\mathfrak{h}$
$$\widetilde{W}\doteq W\ltimes\Lambda(\Phi).$$
$\widetilde{W}$ is called the \emph{extended affine Weyl group}
of $\Phi$ and contains the {affine Weyl group}
$\widehat{W}\doteq W\ltimes\langle\Phi^\vee\rangle$ (see for
instance \cite{IM}, \cite{Ra}).

The action of $W\times Z(\Phi)$ on $\PT$ is induced by that of
$\widetilde{W}$. Indeed $\widetilde{W}$ preserves the lattice
$\langle\Phi^\vee\rangle$ of $\mathfrak{h}$, and thus acts on
$T=\Tquo$ and on $\PT\subset T$ . Since the semidirect factor
$\langle\Phi^\vee\rangle$ acts trivially, $\widetilde{W}$ acts
as its quotient
$$\frac{\widetilde{W}}{\langle\Phi^\vee\rangle}\simeq W\times Z(\Phi).$$\\
\end{re}

\subsection{Examples: the classical root systems}
In the following examples we denote by $\sym_n$,
$\mathfrak{D}_n$, $\mathfrak{C}_n$ respectively the symmetric,
dihedral and cyclic group on $n$ letters.

\begin{enumerate}
\item \textbf{Case $\mathsf{C}_n$} The roots
    $$2\alpha_i+\dots+2\alpha_{n-1}+\alpha_n$$
    $(i=1,\dots,n)$ take integer values on the points
    $[{\alpha^\vee_1}/2],\dots,[{\alpha^\vee_n}/2]\in\Tquo$,
    and thus on their sums, for a total of $2^n$ points of
    $\PT$. Indeed, let us introduce the following notation.
    Fixed a basis $h_1^*,\dots,h_n^*$ of $\mathfrak{h}^*$,
    the simple roots of $\mathsf{C}_n$ can be written as
    \begin{equation}
    \alpha_i=h_i^*-h_{i+1}^* \mbox{ for } i=1,\dots,n-1 \; ,\mbox{ and }\alpha_n=2h_n^*.
    \end{equation}
     Then
     $$\Phi=\{h_i^*-h_j^*\}\cup\{h_i^*+h_j^*\}\cup\{\pm
2h_i^*\}\;(i,j=1,\dots,n\:,\:i\neq j)$$
 and writing $t_i$
for $e^{h_i^*}$, we have that
$$e^\Phi\doteq\{e^\alpha,\alpha\in\Phi\}=\{t_i
t_j^{-1}\}\cup\{t_i t_j\}\cup\{t_i ^{\pm 2}\}.$$ The system
of $n$ independent equations
$$\begin{cases}
t_1^2=1\\
\dots\\
t_n^2=1
\end{cases}$$
has $2^n$ solutions:
 $(\pm1,\dots,\pm1)$, and it is easy to see that all other
systems does not have other solutions. The Weyl group
$W\simeq\sym_n\ltimes (\mathfrak{C}_2)^n$ acts on
$T={(\mathbb{C}^*)}^n$ by permuting and inverting its
coordinates; the second operation is trivial on $\PT$. Thus
two elements of $\PT$ are in the same $W-$orbit if and only
if they have the same number of negative coordinates. Then
we can define the $p-$th $W-$orbit $\mathcal{O}_p$ as the
set of points with $p$ negative coordinates. (This choice
is not canonical: we may choose the set of points with $p$
positive coordinates as well). Clearly if $t\in
\mathcal{O}_p$ then
$$W(t)\simeq(\sym_p\times\sym_{n-p})\ltimes
(\mathfrak{C}_2)^n.$$ Thus $|\mathcal{O}_p|={n \choose p}$
and we get:
$$|\PT|=\sum_{p=0}^n {n \choose p}=2^n.$$
Notice that if $t\in \mathcal{O}_p$ then
$-t\in\mathcal{O}_{n-p}$, and $Ad(t)=Ad(-t)$ since
$Z(\Phi)=\{ \pm(1,\dots,1)\}$. In fact $\Gamma$ has a
symmetry exchanging the vertices $p$ and $n-p$. Finally
notice that $\PT$ is a subgroup of $T$ isomorphic to
${(\mathfrak{C}_2)}^n$ and generated by the elements
$$\delta_i\doteq(1,\dots,1,-1,1,\dots,1) \mbox{ (with the
$-1$ at the $i-th$ place)}.$$ Then we can come back to the
original coordinates observing that $\delta_i$ is the
nontrivial solution of the system ${t_i}^2=1$,
$t_j=1\forall j\neq i$, and using (6) to get:
$$\delta_i\leftrightarrow\left[\sum_{k=i}^n
\alpha_k^\vee/2\right].$$

\item \textbf{Case $\mathsf{D}_n$} We can write
    $\alpha_n=h_{n-1}^*+h_{n}^*$ and the others $\alpha_i$
    as before; then $$e^\Phi=\{t_i t_j^{-1}\}\cup\{t_i
    t_j\}.$$ Then each system of $n$ independent equations
    is $W-$conjugate to one of this form:
    $$\begin{cases}
t_1=t_2\\
\dots\\
t_{p-1}=t_p\\
t_{p-1}=t_p^{-1}\\
t_{p+1}^{\pm 1}=t_{p+2}\\
\dots\\
t_{n-1}=t_n\\
t_{n-1}=t_n^{-1}
\end{cases}$$
for some $p\neq 1,n-1$. Then we get the subset of
${(\mathfrak{C}_2)}^n$ composed by the following $n-$ples:
$$\{(\pm1,\dots,\pm1)\}\setminus
\{\pm\delta_i,\: i=1,\dots,n\}$$ which are in number of $2^n-2n$. However
reasoning as before we see that each one represents two
points in $\Tquo$. Namely, the correspondence is given by:
$$\left\{\left[\sum_{k=i}^{n-1} \frac{\alpha_k^\vee}{2}\pm\frac{\alpha_{n-1}^\vee-\alpha_n^\vee}{4} \right]\right\}\longrightarrow\delta_i.$$
From a geometric point of view, the $t_i$s are coordinates
of a maximal torus of the orthogonal group, while $T=\Tquo$
is a maximal torus of its two-sheets universal covering.
Each $W-$orbit corresponding to the four extremal vertices
of $\Gamma$ is a singleton consisting of one of the four
points over $\pm(1,\dots,1)$, all inducing the identity
automorphism: indeed $Aut(\Gamma)$ acts transitively on
these points. The other orbits are defined as in the case
$\mathsf{C}_n$.

\item \textbf{Case $\mathsf{B}_n$} This case is very
    similar to the previous one, but now $\alpha_n=~h_n^*$,
    $$e^\Phi=\{t_i t_j^{-1}\}\cup\{t_i t_j\}\cup\{t_i^{\pm
    1}\}$$ and then we get the points
    $$\{(\pm1,\dots,\pm1)\}\backslash
    \{\delta_i\}_{i=1,\dots,n}.$$ In this case the projection is
    $$\left\{\left[\sum_{k=i}^{n-1}
    \frac{\alpha_k^\vee}{2}\pm\frac{\alpha_n^\vee}{4}
    \right]\right\}\longrightarrow\delta_i$$ then we have
    $2^n-n$ pairs of points in $\PT$.

\item \textbf{Case $\mathsf{A}_n$} If we see
    $\mathfrak{h}^*$ as the subspace of $\langle
    h_1^*,\dots,h_{n+1}^*\rangle $ of equation $\sum
    h_i^*=0$, and $T$ as the subgroup of
    $(\mathbb{C}^*)^{n+1}$ of equation $\prod t_i=~1$, we
    can write all the simple roots  as
    $\alpha_i=h_i^*-h_{i+1}^*$; then $e^\Phi=~\{t_i
    t_j^{-1}\}$. In this case $\Phi$ has no proper
    subsystem of its same rank, then all the coordinates
    must be equal. Therefore
$$\PT=Z(\Phi)=\left\{(\zeta,\dots,\zeta)|\zeta^{n+1}=1\right\}\simeq \mathfrak{C}_{n+1}.$$
Then $W\simeq\sym_{n+1}$ acts on $\PT$ trivially and
$Z(\Phi)$ transitively, as expected since
$Aut(\Gamma)\simeq \mathfrak{D}_{n+1}$ acts transitively on
the vertices of $\Gamma$. We can write more explicitly
$\PT\subseteq\Tquo$ as
$$\PT=\left\{\left[\frac{k}{n+1}\sum_{i=1}^n i\alpha_i^\vee\right], k=0,\dots,n \right\}.$$\\

\end{enumerate}

\subsection{Proofs}

Motivated by Remark \ref{re2}, we start to describe the
automorphisms of $\mathfrak{g}$ that are induced by the points
of $\PT$.

\begin{lem}
If $t\in\PT$, then $Ad(t)$ has finite order.
\end{lem}
\begin{proof}
Let $\beta_1,\dots,\beta_n$ linearly independent roots such
that $e^{\beta_i}(t)=1$: then for each root $\alpha\in\Phi$ we
have that $m\alpha=\sum c_i\beta_i$ for some $m$ and
$c_i\in\mathbb{Z}$, and thus
$$e^{\alpha}(t^m)=e^{m\alpha}(t)=\prod_{i=1}^n (e^{\beta_i})^{c_i}(t)=1.$$
Then $Ad(t^m)$ is the identity on $\mathfrak{g}$, hence by
(\ref{adg}) $t^m\in Z(\Phi)$. $Z(\Phi)$ is a finite group, thus
$t^m$ and $t$ have finite order.
\end{proof}

The previous lemma allows us to apply the following
\begin{te}[Ka\v{c}]\label{kac}~
\begin{enumerate}
\item Each inner automorphism of $\mathfrak{g}$ of finite
    order $m$  is conjugate to an automorphism $\sigma$ of
    the form
$$\sigma(X_i)=\zeta^{s_i}X_i$$
with $\zeta$ fixed primitive $m-$th root of unity and
$(s_0,\dots,s_n)$ nonnegative integers without common
factors such that $m=\sum s_i a_i$.
\item Two such automorphisms are conjugate if and only if
    there is an automorphism of $\Gamma$ sending the
    parameters $(s_0,\dots,s_n)$ of the first in the
    parameters $(s'_0,\dots,s'_n)$ of the second.
\item Let $(i_1,\dots,i_r)$ be all the indices for which
    $s_{i_1}=\dots=s_{i_r}=0$. Then $\mathfrak{g}^\sigma$
    is the direct sum of an $(n-r)-$dimensional center and
    of a semisimple Lie algebra whose Dynkin diagram is the
    subdiagram of $\Gamma$ of vertices $i_1,\dots,i_r$.
\end{enumerate}
\end{te}
This is a special case of a theorem proved in \cite{KA} and
more extensively in \cite[X.5.15 and 16]{He}. We only need the
following

\begin{co}\label{kco}~
\begin{enumerate}
\item Let $\sigma$ be an inner automorphism of
    $\mathfrak{g}$ of finite order $m$ such that
    $\mathfrak{g}^\sigma$ is semisimple. Then there is
    $p\in V(\Gamma)$ such that $\sigma$ is conjugate to
    $\sigma_p$. In particular $m=a_p$ and the Dynkin
    diagram of $\mathfrak{g}^\sigma$ is $\Gamma_p$.
\item Two automorphisms $\sigma_p$, $\sigma_{p'}$ are
    conjugate if and only if $p,p'$ are in the same
    $Aut(\Gamma)-$orbit.
\end{enumerate}
\end{co}
\begin{proof}
If $\mathfrak{g}^\sigma$ is semisimple, then in the third part
of Theorem \ref{kac} $n=r$, hence all  parameters of $\sigma$
but one are equal to 0, and the nonzero parameter $s_p$ must be
equal to 1, otherwise there would be a common factor,
contradicting the first part of the Theorem. Thus we get the
first statement. Then the second statement follows from Theorem
\ref{kac}(ii).
\end{proof}

Let be $t\in\PT$: by Remark \ref{re2} $\mathfrak{g}^{Ad(t)}$ is
semisimple, hence by Corollary \ref{kco}(i) $Ad(t)$ is
conjugate to some $\sigma_p$. Then there is a canonical map
\begin{align*}
\psi : \PT & \rightarrow Q\\
t & \mapsto \psi(t)=\{ p\in V(\Gamma) \mbox{ such that }\sigma_p \mbox { is conjugate to }Ad(t)\}.
\end{align*}
Notice that $\psi(t)$ is a well-defined element of $Q$ by
Corollary \ref{kco}(ii).

We now prove the fundamental
\begin{lem}
Two points in $\PT$ induce conjugate automorphisms if and only
if they are in the same $W\times Z(\Phi)-$orbit.
\end{lem}

\begin{proof}
We recall that $W\simeq N/T$ and the action of $W$ on $T$ is
induced by the conjugation action of $N$; it is also well known
that two points of $T$ are $G-$conjugate if and only if they
are $W-$conjugate. Then $W-$conjugate points induce conjugate
automorphisms. Moreover by (\ref{adg})
$$Ad(t)=Ad(s)\Leftrightarrow
Ad(ts^{-1})=id_{\mathfrak{g}}\Leftrightarrow ts^{-1}\in
Z(\Phi).$$ Finally suppose that $t,t'\in\PT$ induce conjugate
automorphisms, i.e.
$$\exists g\in G | Ad(t')=Ad(g)Ad(t)Ad(g^{-1})=Ad(gtg^{-1}).$$
Then $zt'=gtg^{-1}$ for some $z\in Z(\Phi)$. Thus $zt'$ and $t$ are $G-$conjugate elements of $T$, and hence they are $W-$conjugate, proving the claim.\\
\end{proof}

We can now prove the first part of Theorem \ref{pca}. Indeed by
the previous lemma there is a canonical injective map defined
on the set of the orbits of $\PT$:
$$\overline{\psi} : \cfrac{\PT}{W\times Z(\Phi)} \longrightarrow Q.$$
We must show that this map is surjective. The system
$$\alpha_i(h)=1 \,(\forall i\neq 0,p)\, \mbox{, } \alpha_p(h)=a_p^{-1}$$
 is composed of $n$ linearly independent equations, then it has a solution $h\in\mathfrak{h}$. Notice that $\alpha_0(h)\in\mathbb{Z}$.
 Let $t$ be the class of $h$ in $T$; then $$e^\alpha(t)=1\Leftrightarrow\alpha\in \Phi_p.$$
  Then by Remark \ref{re2} $Ad(t)$ is conjugate to $\sigma_p$ and $\Phi(t)$ to $\Phi_p$.\\

In order to relate the action of $Z(\Phi)$ with that of
$Aut(\Gamma)$, we introduce the following subset of $W$. For
each $p\neq 0$ such that $a_p=1$, set $z_p\doteq w_0^p w_0$,
where $w_0$ is the longest element of $W$ and $w_0^p$ is the
longest element of the parabolic subgroup of $W$ generated by
all the simple reflections $s_{\alpha_1},\dots,s_{\alpha_n}$
except $s_{\alpha_p}$. Then we define
$$W_Z\doteq \{1\}\cup \{z_p\}_{p=1,\dots,n | a_p=1}$$

$W_Z$ has the following properties (see \cite[§1.7 and
1.8]{IM}):
\begin{te}[Iwahori-Matsumoto]~
\begin{enumerate}
  \item $W_Z$ is a subgroup of $W$ isomorphic to $Z(\Phi)$.
  \item For each $z_p\in W_Z$, we have that
      $z_p.\alpha_0=\alpha_p$, and $z_p$ induces an
      automorphism of $\Gamma$ that sends the $0-$th vertex
      to the $p-$th one; this defines an injective morphism
      $W_Z\hookrightarrow Aut(\Gamma)$.
      \item The $W_Z-$orbits of $V(\Gamma)$ coincide with
          the $Aut(\Gamma)-$orbits.
\end{enumerate}
\end{te}

Therefore $Q$ is the set of $W_Z-$orbits of $V(\Gamma)$, and
the bijection $\overline{\psi}$ between $Q$ and the set of
$Z(\Phi)-$orbits of $\PT/W$ can be lifted to a noncanonical
bijection  between $V(\Gamma)$ and $\PT/W$. Then we just have
to consider the action of $W$ on $\PT$ and prove the
\begin{lem}
If $t\in\mathcal{O}_p$, then $W(t)$ is conjugate to $W_p$.
\end{lem}

\begin{proof}
Notice that the centralizer $C_N(t)$ of $t$ in $N$ is the normalizer of $T=~C_T(t)$ in $C_G(t)$. Then $W(t)=C_N(t)/T$ is the Weyl group of $C_G(t)$. $C_G(t)$ is the subgroup of $G$ of points fixed by the conjugacy by $t$, then its Lie algebra is $\mathfrak{g}^{Ad(t)}$, which is conjugate to $\mathfrak{g}^{\sigma_p}$ by the first part of Theorem \ref{pca}. Therefore $W(t)$ is conjugate to $W_p$.\\
\end{proof}

This completes the proof of Theorem \ref{pca} and also of
Theorem \ref{pts}, since by Remark \ref{re2} the map $\psi$
defined in (7) can also be seen as the map
$$t \mapsto \psi(t)=\{ p\in V(\Gamma) \mbox{ such that } \Phi_p \mbox { is conjugate to }\Phi(t)\}.$$\\

\section{Layers of the arrangement}

\subsection{From hyperplane arrangements to toric arrangements}

 Let $S$ be a $d-$dimensional space of $\mathcal{H}$.
 The set $\Phi_S$ of the elements of $\Phi$ vanishing on $S$ is a complete subsystem of $\Phi$ of rank $n-d$.
 Then the map $S\rightarrow \Phi_S$ gives a bijection between $\mathcal{L}_d$ and $\mathcal{K}_d$, whose inverse is
$$\Theta\rightarrow S(\Theta)\doteq\{h\in \mathfrak{h} | \alpha(h)=0 \; \forall\alpha\in\Theta\}.$$

In \cite[6.4 and C]{OT} (following \cite{OS} and \cite{Ca}) the
spaces of $\mathcal{H}$ are classified and counted, and the
$W-$orbits of $\mathcal{L}_d$ are completely described. This is
done case-by-case according to the type of $\Phi$. We now show
a case-free way to extend this analysis to the layers of
$\mathcal{T}$.

Given a layer $C$ of $\mathcal{T}$ let us consider
$$\Phi_C\doteq \{\alpha\in\Phi | e^{\alpha}(t)=1 \; \forall
t\in C\}.$$ In contrast with the case of linear arrangements,
$\Phi_C$ in general is not complete. For each $\Theta\in
\mathcal{K}_d$, define $\mathcal{C}^{\Phi}_{\Theta}$ as the set
of layers $C$ such that $\overline{\Phi_C}=\Theta$. This is
clearly a partition of the set of $d-$dimensional layers of
$\mathcal{T}$, i.e.
\begin{equation}\label{pa8}
\mathcal{C}_{d}(\Phi)=\bigsqcup_{\Theta\in \mathcal{K}_d} \mathcal{C}^{\Phi}_{\Theta}
\end{equation}
Given any $C\in\mathcal{C}^{\Phi}_{\Theta}$, we call
$S(\Theta)$ the \emph{tangent space} at the layer $U$.
Then by \cite{OT} the problem of classifying the layers of $\mathcal{T}$ reduces to classify the layers of $\mathcal{T}$
having a given tangent space, i.e. the elements of $\mathcal{C}^{\Phi}_{\Theta}$.
In the next section we show that this amounts to classify the points of a smaller toric arrangement, namely that defined by $\Theta$.\\

\subsection{Theorems}

Let $\Theta$ be a complete subsystem of $\Phi$ and $W^{\Theta}$
its Weyl group. Let $\mathfrak{k}$ and $K$ be respectively the
semisimple Lie algebra and the semisimple and simply connected
algebraic group of root system $\Theta$, $\mathfrak{d}$ a
Cartan subalgebra of $\mathfrak{k}$,
$\langle\Theta^\vee\rangle$ and $\Lambda(\Theta)$ the coroot
and coweight lattices, $Z(\Theta)\doteq~
\frac{\Lambda(\Theta)}{\langle\Theta^\vee\rangle}$ the center
of $K$, $D$ the maximal torus of $K$ defined by
${\mathfrak{d}}/{\langle\Theta^\vee\rangle}$, $\mathcal{D}$ the
toric arrangement defined by $\Theta$ on $D$ and
$\mathcal{C}_0(\Theta)$ the set of its points.

We also consider the \emph{adjoint group} $K_a\doteq
K/Z(\Theta)$ and its maximal torus $D_a\doteq D/Z(\Theta)\simeq
\mathfrak{d}/\Lambda(\Theta)$. We recall from \cite{Hu} that
$K$ is the universal covering of $K_a$, and if $D'$ is an
algebraic torus having Lie algebra $\mathfrak{d}$, then
$D'\simeq \mathfrak{d}/L$ for some lattice
$\Lambda(\Theta)\supseteq L \supseteq
\langle\Theta^\vee\rangle$; then there are natural covering
projections $D\twoheadrightarrow D'\twoheadrightarrow D_a$ with
kernels respectively $L/\langle\Theta^\vee\rangle$ and
$\Lambda(\Theta)/L$. Notice that $\Theta$ naturally defines an
arrangement on each torus $D'$, and that for $D'=D_a$ the set
of its 0-dimensional layers is
$\mathcal{C}_{0}(\Theta)/Z(\Theta)$. Given a point $t$ of some
$D'$ we set
$$\Theta(t)\doteq \{\alpha\in\Theta|e^\alpha(t)=1\}.$$

\begin{te}\label{cft}
There is a $W^\Theta-$equivariant surjective map
$$\varphi: \mathcal{C}^{\Phi}_{\Theta} \twoheadrightarrow  \mathcal{C}_{0}(\Theta)/Z(\Theta)$$
such that $\ker \varphi\simeq Z(\Phi)\cap Z(\Theta)$ and
$\Phi_C=\Theta(\varphi(C))$.
\end{te}

\begin{proof}
Let $S(\Theta)$ be the subspace of $\mathfrak{h}$ defined in
the previous section, and $H$ the corresponding subtorus of
$T$. $T/H$ is a torus with Lie algebra
$\mathfrak{h}/S(\Theta)\simeq~\mathfrak{d}$, then $\Theta$
defines an arrangement $\mathcal{D}'$ on $D'\doteq T/H$. The
projection $\pi:T\twoheadrightarrow T/H$ induces a bijection
between $\mathcal{C}^{\Phi}_{\Theta}$ and the set of
0-dimensional layers of $\mathcal{D}'$, because
$H\in\mathcal{C}^{\Phi}_{\Theta}$ and for each
$C\in\mathcal{C}^{\Phi}_{\Theta}$, $\Phi_C=\Theta(\pi(C))$.

Moreover the restriction of the projection $d\pi:\mathfrak{h}
\twoheadrightarrow \mathfrak{h}/S(\Theta)$ to
$\langle\Phi^\vee\rangle$ is simply the map that restricts the
coroots of $\Phi$ to $\Theta$. Set $R^{\Phi}(\Theta)\doteq
d\pi(\langle\Phi^\vee\rangle)$; then $\Lambda(\Theta)\supseteq
R^{\Phi}(\Theta) \supseteq \langle\Theta^\vee\rangle$ and
$D'\simeq \mathfrak{d}/R^{\Phi}(\Theta)$. Denote by $p$ the
projection
$\Lambda(\Phi)\twoheadrightarrow\frac{\Lambda(\Phi)}{\langle\Phi^\vee\rangle}$
and embed $\Lambda(\Theta)$ in $\Lambda(\Phi)$ in the natural
way. Then the kernel of the covering projection of
$D'\twoheadrightarrow D_a$ is isomorphic to
$$\frac{\Lambda(\Theta)}{R^{\Phi}(\Theta)}\simeq p(\Lambda(\Theta))\simeq Z(\Phi)\cap Z(\Theta).$$
\end{proof}

We set
$$n_\Theta\doteq \frac{|Z(\Theta)|} {|Z(\Phi)\cap Z(\Theta)|}.$$
The following corollary is straightforward from Theorem
\ref{cft}.
\begin{co}\label{nla}
$$|\mathcal{C}^{\Phi}_{\Theta}|=n_\Theta^{-1} |\mathcal{C}_{0}(\Theta)|$$
and then by (\ref{pa8}),
$$|\mathcal{C}_{d}(\Phi)|=\sum_{\Theta\in \mathcal{K}_d} n_\Theta^{-1} |\mathcal{C}_{0}(\Theta)|.$$
\end{co}

Notice that two layers $C, C'$ of $\mathcal{T}$ are
$W-$conjugate if and only if the two conditions below are
satisfied:
\begin{enumerate}
  \item  their tangent spaces are $W-$conjugate , i.e.
      $\exists w\in W$ such that
      $\overline{\Phi_C}=w.\overline{\Phi_{C'}}$;
  \item $C$ and $w.C'$ are
      $W^{\overline{\Phi_C}}-$conjugate.
\end{enumerate}
 Then the action of $W$ on $\mathcal{C}(\Phi)$ is described by the following remark.

\begin{re}~
\begin{enumerate}
  \item By Theorem \ref{cft}, $\varphi$ induces a
      surjective map $\overline{\varphi}$ from the set of
      the $W^\Theta-$orbits of
      $\mathcal{C}^{\Phi}_{\Theta}$ to the set of the
      $W^\Theta\times Z(\Theta)-$orbits of
      $\mathcal{C}_{0}(\Theta)$, that are described by
      Theorem \ref{pca}.

\item In particular if $\Theta$ is irreducible, set
    $\Gamma^{\Theta}$ its affine Dynkin diagram, $Q^\Theta$
    the set of the $Aut(\Gamma)-$orbits of its vertices,
    $\Gamma^{\Theta}_p$ the diagram that we obtain from
    $\Gamma^{\Theta}$ removing the vertex $p$, and
    $\Theta_p$ the associated root system. Then there is a
    surjective map
  $$\widehat{\varphi}:\mathcal{C}^{\Phi}_{\Theta}\twoheadrightarrow Q^\Theta$$
  such that, if $\widehat{\varphi}(C)=q$ and $p$ is a representative of $q$, then $\Phi_C\simeq \Theta_p$.\\
\end{enumerate}
\end{re}

\subsection{Examples}

\textbf{Case $\mathsf{F}_4$. } $Z(\Phi)=\{1\}$, thus
$n_\Theta=|Z(\Theta)|$. Therefore in this case $n_\Theta$ does
not depend on the conjugacy class, but only on the isomorphism
class of $\Theta$.

We say that a space $S$ of $\mathcal{H}$ (respectively a layer
$C$ of $\mathcal{T}$) is of a given type  if the corresponding
subsystem $\Phi_S$ (respectively $\Phi_C$) is of that type.
Then by \cite[Tab. C.9]{OT} and Corollary \ref{nla} there are:
\begin{enumerate}
  \item one space of type $"\mathsf{A_0}"$, tangent to one
      layer of the same type (the whole spaces);
  \item 24 spaces of type $\mathsf{A_1}$, each tangent to
      one layer of the same type;
  \item 72 spaces of type $\mathsf{A_1}\times
      \mathsf{A_1}$, each tangent to one layer of the same
      type;
  \item 32 spaces of type $\mathsf{A_2}$, each tangent to
      one layer of the same type;
  \item 18 spaces of type $\mathsf{B_2}$, each tangent to
      one layer of the same type and one layer of type
      $\mathsf{A_1}\times \mathsf{A_1}$;
  \item 12 spaces of type $\mathsf{C_3}$, each tangent to
      one layer of the same type and 3 of type
      $\mathsf{A_2}\times \mathsf{A_1}$;
  \item 12 spaces of type $\mathsf{B_3}$, each tangent to
      one layer of the same type, one of type
      $\mathsf{A_3}$ and 3 of type $\mathsf{A_1}\times
      \mathsf{A_1}\times \mathsf{A_1}$;
  \item 96 spaces of type $\mathsf{A_1}\times
      \mathsf{A_2}$, each tangent to one layer of the same
      type;
  \item one space of type $\mathsf{F_4}$ (the origin),
      tangent to: one layer of the same type, 12 of type
      $\mathsf{A_1}\times \mathsf{C_3}$,
      32 of type $\mathsf{A_2}\times \mathsf{A_2}$, 24 of
      type $\mathsf{A_3}\times \mathsf{A_1}$,
      and 3 of type $\mathsf{C_4}$.\\
      \end{enumerate}

\textbf{Case $\mathsf{A}_{n-1}$. } It is easily seen that each
subsystem $\Theta$ of $\Phi$ is complete and is a product of
irreducible factors $\Theta_1,\dots,\Theta_k$, with $\Theta_i$
of type $A_{\lambda_i-1}$ for some positive integers
$\lambda_i$ such that $\lambda_1+\dots +\lambda_k=n$  and $n-k$
is the rank of $\Theta$. In other words, as is well known, the
$W-$conjugacy classes of spaces of $\mathcal{H}$ are in
bijection with the partitions $\lambda$ of $n$, and if a space
has dimension $d$ then corresponding partition has length
$|\lambda|\doteq k$ equal to $d+1$. The number of spaces of
partition $\lambda$ is easily seen to be equal to
$n!/b_\lambda$, where $b_i$ is the number of $\lambda_j$ that
are equal to $i$ and $b_\lambda\doteq \prod i!^{b_i}b_i!$
 (see \cite[6.72]{OT}). Now let $g_\lambda$ be the greatest common divisor of $\lambda_1,\dots,\lambda_k$. By Example 4 in Section 2.2
 we have that $$|Z(\Theta)|=\lambda_1 \dots  \lambda_k=|\mathcal{C}_0(\Theta)|$$ and $|Z(\Phi)\cap Z(\Theta)|=g_\lambda$. Then by Corollary \ref{nla} $|\mathcal{C}^{\Phi}_{\Theta}|=g_\lambda$ and $$|\mathcal{C}_d(\Phi)|=\sum_{|\lambda|=d+1}\frac{n! g_\lambda}{b_\lambda}.$$

This could also be seen directly as follows. We can view $T$ as
the subgroup of $(\mathbb{C}^*)^{n}$ given by the equation
$t_1\dots t_n -1=0$. Then $\Theta$ imposes the equations
$$\begin{cases}
t_1=\dots=t_{\lambda_1}\\
\dots\\
t_{\lambda_1+\dots+\lambda_{k-1}+1}=\dots=t_n.
\end{cases}$$
Thus we have the relation $$x_1^{\lambda_1}\dots
x_k^{\lambda_k}-1=0.$$ If $g_\lambda=1$ this polynomial is
irreducible, because the vector $(\lambda_1,\dots,\lambda_k)$
can be completed to a basis of the lattice $\mathbb{Z}^{k}$. If
$g_\lambda>1$ this polynomial has exactly $g_\lambda$
irreducible factors over $\mathbb{C}$. Then in every case it
defines an affine variety having $g_\lambda$ irreducible
components,
which are precisely the elements of $\mathcal{C}^{\Phi}_{\Theta}$.\\

\section{Topology of the complement}
\subsection{Theorems}

Let $\rp$ be the complement of the toric arrangement:
$$\rp\doteq T\setminus\bigcup_{\alpha\in \Phi^+}H_\alpha.$$

In this section we prove that the Euler characteristic of
$\rp$, denoted by $E_\Phi$, is equal to $(-1)^n |W|$. This may
also be seen as a consequence of \cite[Prop. 5.3]{CMS}.
Furthermore, we give a formula for the Poincar\'{e} polynomial
of $\rp$, denoted by $P_\Phi(q)$.

Let $d_1,\dots,d_n$ be the \emph{degrees} of $W$, i.e. the
degrees of the generators of the ring of $W-$invariant regular
functions on $\mathfrak{h}$; it is well known that $d_1 \dots
d_n=|W|$. The numbers $d_1-1,\dots,d_n-1$ are known as the
\emph{exponents} of $W$; we denote by $\mathcal{P}(\Phi)$ their
product:
 $$\mathcal{P}(\Phi)\doteq(d_1-1)\dots(d_n-1).$$
Then we have:

\begin{te}
    $$P_\Phi(q)=\sum_{C\in \mathcal{C}(\Phi)} \mathcal{P}(\Phi_C) (q+1)^{d(C)} q^{n-d(C)}$$
where $d(C)$ is the dimension of the layer $C$.
\end{te}
\begin{proof}
Let $nbc(\Phi)$ be the number of \emph{no-broken circuit bases}
of $\Phi$: by
\cite{OSb}, $nbc(\Phi)$ equals the leading coefficient of the
Poincar\'{e} polynomial of the complement of $\mathcal{H}$ in
$\mathfrak{h}$; moreover by \cite{Br} this coefficient is equal
to $\mathcal{P}(\Phi)$ (these facts can be found also in
\cite[10.1]{li}).

 Then the claim is a restatement of a known result. Indeed the cohomology of $\rp$ can be expressed as a
direct sum of contributions given by the layers of
$\mathcal{T}$ (see for example \cite[Theor. 4.2]{DPt} or
\cite[14.1.5]{li}). In terms of Poincar\'{e} polynomial this
expression is:
$$P_\Phi(q)=\sum_{C\in \mathcal{C}(\Phi)} nbc(\Phi_C) (q+1)^{d(C)} q^{n-d(C)}.$$
\end{proof}

Now we use the theorem above to compute the Euler
characteristic of $\rp$.
\begin{lem}
$$E_\Phi= (-1)^n\sum_{p=0}^n \frac{|W|}{|{W_p}|}\mathcal{P}(\Phi_p)$$
\end{lem}

\begin{proof}
We have
\begin{equation}
E_\Phi=P_\Phi(-1)=(-1)^n\sum_{t\in\PT}\mathcal{P}(\Phi(t))
\end{equation}
 because the contributions of all positive-dimensional layers vanish at $-1$.
 Obviously isomorphic subsystems have the same degrees, thus Theorem \ref{pts} yields the statement.
\end{proof}

\begin{te}\label{euf}
$$E_\Phi=(-1)^n |W|$$
\end{te}

\begin{proof}
By the previous lemma we must prove that $$\sum_{p=0}^n
\frac{\mathcal{P}(\Phi_p)}{|{W_p}|}=1$$ If we write
$d_1^p,\dots,d_n^p$ for the degrees of ${W_p}$, the previous
identity becomes
$$\sum_{p=0}^n \frac{(d_1^p-1)\dots(d_n^p-1)}{d_1^p\dots d_n^p}=1.$$
This identity has been proved in \cite{ci}, and later with
different methods in \cite{De}.
\end{proof}

Notice that $W$ acts on $\rp$ and then on its cohomology. Then
we can consider the \emph{equivariant Euler characteristic} of
$\rp$, that is, for each $w\in W$,
$$\widetilde{E}_\Phi(w)\doteq \sum_{i=0}^{n} (-1)^i \: Tr(w,H^{i}(\rp,\mathbb{C}) ).$$
Let $\varrho_W$ be the character of the regular representation
of $W$. From Theorem \ref{euf} we get the following

\begin{co}
$$\widetilde{E}_\Phi=(-1)^n \varrho_W$$
\end{co}

\begin{proof}
Since $W$ is finite and acts freely on $\rp$, it is well known
that $\widetilde{E}_\Phi=k \varrho_W$ for some $k\in
\mathbb{Z}$. Then to compute $k$ we just have to look at
$\widetilde{E}_\Phi(1_W)=~E_\Phi$.
\end{proof}

Finally we give a formula for $P_\Phi(q)$ which, together with
the mentioned results in \cite{OT}, allows its explicit
computation.

\begin{te}
$$P_\Phi(q)=\sum_{d=0}^n (q+1)^d q^{n-d }\sum_{\Theta\in \mathcal{K}_d}n_\Theta^{-1}|W^\Theta|$$
\end{te}

\begin{proof}
By formula (\ref{pa8}) we can restate Theorem 14 as
$$P_\Phi(q)=\sum_{d=0}^n (q+1)^d q^{n-d }\sum_{\Theta\in \mathcal{K}_d}\sum_{C\in \mathcal{C}^\Phi_\Theta} \mathcal{P}(\Phi_C)$$
Moreover by Theorem \ref{cft} and Corollary \ref{nla} we get
$$\sum_{C\in \mathcal{C}^\Phi_\Theta} \mathcal{P}(\Phi_C)=n_\Theta^{-1}\sum_{t\in\Pt}\mathcal{P}(\Theta(t)).$$
Finally the claim follows by formula (9) and Theorem
\ref{euf} applied to $\Theta$:
$$\sum_{t\in\Pt}\mathcal{P}(\Theta(t))=(-1)^d \chi_\Theta=|W^\Theta|.$$\\
\end{proof}

\subsection{Examples}

\textbf{Case $\mathsf{F}_{4}$. } In Section 3.3 we have given
a list of all possible types of complete subsystems, together
with their multiplicities. Then we just have to compute the
coefficient $n_\Theta^{-1}|W^\Theta|$ for each type. This is
equal to:
\begin{itemize}
  \item 1 for types 1., 2. and 3.
  \item 2 for types 4. and 8.
  \item 4 for type 5.
  \item 24 for types 6. and 7.
  \item 1152 for type 9.
\end{itemize}

Thus
$$P_\Phi(q)=2153 q^4 +1260 q^3+ 286 q^2 +28 q+1.$$\\

\textbf{Case $\mathsf{A}_{n-1}$. } By Section 1.3.3,
$n_\Theta^{-1}=\frac{g_\lambda}{\lambda_1\dots\lambda_k}$ and
$|W^\Theta|=\lambda_1!\dots \lambda_k!$. Hence by Theorem 17
$$P_\Phi(q)=\sum_{d=0}^n (q+1)^d q^{n-d }\sum_{|\lambda|=d+1}n!b_\lambda^{-1}g_\lambda(\lambda_1-1)!\dots (\lambda_k-1)!.$$\\

\end{document}